\documentclass[12pt]{article}
\usepackage{latexsym,amssymb}
\usepackage{latexsym}
\usepackage{graphicx}
\usepackage[labelfont=bf]{caption}
\newtheorem{thm}{Theorem}[section]

\newtheorem{exs}{Examples}[section]

\newtheorem{rem}{Remark}[section]

\newfont{\sss}{cmss12 scaled 1000}

\title{Simple $t$-designs: A recursive construction for arbitrary $t$}
\author{ Tran van Trung  \\
         Institute for Experimental Mathematics \\ 
         University of Duisburg-Essen \\
         Thea-Leymann-Stra\ss e 9, 45127 Essen, Germany }

\date{}

\begin{document}                  
                
\maketitle

\begin{abstract}
  The aim of this paper is to present a recursive construction of 
  simple $t-$designs for arbitrary $t$. 
  The construction is of purely combinatorial nature
  and it requires finding solutions for the indices of the ingredient 
  designs that satisfy a certain set of equalities. 
  We give a small number of examples to illustrate the construction, whereby 
  we have found a large number of new $t-$designs, which were previously 
  unknown. 
  This indicates that the method is useful and powerful.
  
 \end{abstract}

\vspace{0.1in}\noindent

\vspace{0.1in}\noindent


\section{Introduction}
One of the most challenging problems in design theory is  
the problem of constructing simple $t-$designs for large $t$. 
There are several major approaches to the problem. 
 These are 
constructing $t-$designs from large sets of $t-$designs, for instance
\cite{Ajoodani96}, \cite{Khos-Ajoodani91},
 \cite{Kramer93}, \cite{Kreher93}, 
\cite{Teirlinck87}, \cite{Teirlinck89}, \cite{QRWu91};
constructing $t-$designs by using prescribed automorphism groups,
for example
\cite{betten95}, \cite{betten98}, \cite{betten99},
 \cite{bier95a}, \cite{bier95b},
\cite{Denniston76}, \cite{Kramer76}, \cite{Mag84}; 
or contructing $t-$designs via recursive construction methods,
 see for instance
\cite{driessen78}, \cite{Jimbo2011}, \cite{Mag87}, \cite{Tuan2002}, 
\cite{Sebille2001}, \cite{TvT84}, \cite{TvT86}, \cite{TvT2001}.

In this paper we present a new recursive method for constructing 
simple $t-$designs for arbitrary $t$. The method
is of combinatorial nature, which is a composition technique where
a $t-$design is built up from other smaller ingredient designs.
Which ingredient designs will be necessary are determined by the
solutions to a set of equalities involving their indices. The method
proves to be very useful and powerful.  
Our experimental results obtained from its application 
have shown that, even for a small number of chosen parameters
for the ingredient designs, 
plentiful new simple designs can be constructed, which were previously
unknown.

We recall some basic definitions. A $t-$design, denoted by
$t-(v,k,\lambda)$, is a pair $(X, {\cal B})$, where $X$ is a 
$v-$set of {\it points} and ${\cal B}$ is a collection of $k-$subsets,
called {\it blocks}, of $X$ having the property that every $t-$set 
of $X$ is a subset of exactly $\lambda$ blocks in ${\cal B}$.
The parameter $\lambda$ is called the {\it index} of the design. 
A $t-$design is called {\it simple} if no two blocks are
 identical
i.e. no block of ${\cal B}$ is repeated; otherwise, it is called
{\it non-simple} (i.e. ${\cal B}$ is a multiset). 
It can be shown by simple counting that a
$t-(v,k,\lambda)$ design is an $s-(v,k,\lambda_s)$ design for 
$0 \leq s \leq t$,
where $\lambda_s = \lambda{v-s \choose t-s}/{k-s \choose t-s}.$
Since $\lambda_s$ is an integer, necessary 
conditions for the parameters of a $t-$design are
${k-s \choose t-s} | \lambda{v-s \choose t-s}$, for $0 \leq s \leq t.$
For given $t,v$ and $k$, we denote by 
$\lambda_{\mbox{min}}(t,k,v)$, or 
$\lambda_{\mbox{min}}$ for short, the smallest
positive integer such that these conditions are satisfied for
all $0 \leq s \leq t.$ By complementing each block in $X$ of a
$t-(v,k,\lambda)$ design, we obtain a $t-(v,v-k, \lambda^*)$ design,
where $ \lambda^* = \lambda{v-k \choose t}/{k \choose t}$, hence we shall 
assume that $k \leq v/2.$ The largest value for $\lambda$ for which
a simple $t-(v,k, \lambda)$ design exists is denoted 
by $\lambda_{\mbox{max}}$
and we have $\lambda_{\mbox{max}}={v-t \choose k-t}.$ The simple
$t-(v,k,\lambda_{\mbox{max}})$ design is called the {\it complete}
design or the {\it trivial} design. A $t-(v,k,1)$ design is called a
{\em t-Steiner system.}

We refer the reader to \cite{bjl99}, 
\cite{handbook07} for more information about designs.


\subsection{The Construction}
  
 We first introduce ingredients and notation used in the construction.

 Let $t, \; v, \; k$ be non-negative integers such that 
 $v \geq k \geq t \geq 0$.
 Let $X$ be a $v$-set and let $X= X_1 \cup X_2$ be a partition of $X$
 (i.e $X_1 \cap X_2 = \emptyset$) with $|X_1|=v_1$ and $|X_2|=v_2$.  
 
 Throughout the paper the parameter set
 $t-(v_2,j, \bar{\lambda}^{(j)}_t)$ for a design indicates that the point set 
 of the design is $X_2$. Also, a design defined on the point set $X_2$ 
 will be denoted by $\bar{D}=(X_2, \bar{\cal{B}})$.

 \begin{itemize}
 \item[1.]
  For $i=0, \ldots, t$, let $D_i =(X_1, {\cal B}^{(i)})$ be the 
  complete $i-(v_1,i,1)$ design. For $i=t+1, \ldots, k$, let 
  $D_i =(X_1, {\cal B}^{(i)})$ be a simple $t-(v_1,i, \lambda^{(i)}_t)$
  design.
 \item[2.]
  Similarly, for $i=0, \ldots, t$, let 
  $\bar{D}_i =(X_2, \bar{\cal B}^{(i)})$ 
  be the  complete $i-(v_2,i,1)$ design. And for $i=t+1, \ldots, k$, 
  let   $\bar{D}_i =(X_2, \bar{\cal B}^{(i)})$ be a 
  simple $t-(v_2,i, \bar{\lambda}^{(i)}_t)$ design.
 \item[3.]
  Two degenerate cases for designs occur when either $v=k=t=0$ or $v=k$. 
  The first case $v=k=t=0$ gives an ``empty''
  design, denoted by $\emptyset$, however we use the convention that the
  number of blocks of the empty design is 1 (i.e. the unique block is
  the empty block). The second case $v=k$ gives a degenerate $k$-design
  having just 1 block consisting of all $v$ points. 
  Thus, in these two extreme
  cases the number of blocks of the designs is always 1.

 \item[4.]
  We denote by $T_{(s,t-s)}$ a $t$-subset $T$ of $X$ with
  $|T\cap X_1|=s$ and hence $|T\cap X_2|=t-s$, for $s=0,\ldots, t$.
  It is clear that any $t$-subset of $X$ is a $T_{(s,t-s)}$ set 
  for some $s \in \{0, \ldots, t\}$.
 
  \item[5.]
  Let $X$ be a finite set and let $u \in\{0,1\}$. The notation $X \times [u]$ 
  has the following meaning. $X \times [0]$ is  
  the empty set $\emptyset$, and  $X \times [1] = X$.

 \end{itemize}
      
 We now describe our construction.
 Consider $(k+1)$ pairs of simple designs $(D_i,\bar{D}_{k-i})$
 for $i=0, \ldots, k$, where   
  $D_i =(X_1, {\cal B}^{(i)})$ is a simple $t-(v_1,i, \lambda^{(i)}_t)$
  design and 
  $\bar{D}_{k-i} =(X_2, \bar{\cal B}^{(k-i)})$
  a simple $t-(v_2,k-i, \bar{\lambda}^{(k-i)}_t)$ design,
  as defined above. 
 For each pair $(D_i,\bar{D}_{k-i})$ define
  $${\cal B}_{(i,{k-i})} := \{ B= B_i \cup \bar{B}_{k-i} \; / \;
   B_i \in {\cal B}^{(i)}, \bar{B}_{k-i} \in \bar{\cal B}^{(k-i)} \}. $$ 
 Thus, ${\cal B}_{(i,{k-i})}$ is a collection of $k$-subsets of $X$
 obtained by taking the union of blocks of
 $D_i$ and $\bar{D}_{k-i}$. Note that
 the sets ${\cal B}_{(i,{k-i})}$ and ${\cal B}_{(j,{k-j})}$
 are pairwise disjoint for $i \not= j$ and $i,j = 0, \ldots, k.$ 

 Define
 $${\cal B} := {\cal B}_{(0,k)}\times [u_0] \cup 
             {\cal B}_{(1,k-1)}\times [u_1] 
             \cup   \cdots \cup
             {\cal B}_{(k-1,1)} \times [u_{k-1}] 
             \cup {\cal B}_{(k,0)}\times [u_k], $$ 
where $u_i \in \{0,1\}$, for $i=0, \ldots, k.$  

 It should be noted that the notation 
 ${\cal B}_{(i,k-i)}\times [u_i]$, as defined in [5.], 
 indicates that either we
 have an empty set $\emptyset$ (when $u_i=0$) or the set 
 ${\cal B}_{(i,k-i)}$ itself (when $u_i=1$).
 The empty set case implies that the pair $(D_i,\bar{D}_{k-i})$
 is not used and the other case shows the use of $(D_i,\bar{D}_{k-i})$.
 Thus $u_i$'s are considered as variables.

 We examine the necessary conditions for which $(X, {\cal B})$ forms
 a simple $t$-design. Consider the block set ${\cal B}_{(i,k-i)}.$ 
 We see that
 each $t$-subset $T_{(s,t-s)}$ of $X$ is contained in
  $$    \lambda^{(i)}_s. \bar{\lambda}^{(k-i)}_{t-s}  $$
 blocks of ${\cal B}_{(i,k-i)}$, for $s=0, \ldots, t$. 
 It is clear because any $s$-set of $X_1$ is contained in 
 $\lambda^{(i)}_s$ blocks
 of $D_i$ and any $(t-s)$-set of $X_2$ is contained in  
 $\bar{\lambda}^{(k-i)}_{t-s}$ blocks of $\bar{D}_{k-i}$.
 Note that $\lambda^{(i)}_s. \bar{\lambda}^{(k-i)}_{t-s}$
 could be equal to 0; this is the case when $ i < s$ or $k-i < t-s.$
 Define 
 $$\Lambda^{(i,k-i)}_{s,t-s}:= \lambda^{(i)}_s. \bar{\lambda}^{(k-i)}_{t-s}.$$
 
 It follows that for a given $t$-set $T_{(s,t-s)}$ of $X$ the number 
 of blocks in $\cal B$ containing $T_{(s,t-s)}$ is equal to

\begin{eqnarray*}
 L_{s,t-s} & := & u_0. \Lambda^{(0,k)}_{s,t-s}
               + u_1.\Lambda^{(1,k-1)}_{s,t-s} 
               + \cdots + u_k.\Lambda^{(k,0)}_{s,t-s}  \nonumber \\
           &=& \sum_{i=0}^k u_i. \Lambda^{(i,k-i)}_{s,t-s} \nonumber \\
           &=& \sum_{i=0}^k u_i.\lambda^{(i)}_s. \bar{\lambda}^{(k-i)}_{t-s}, 
\end{eqnarray*}

Since any $t$-set $T$ of $X$ is of form $T_{s,t-s}$ for some
$s \in \{0, \ldots, t\}$, so if
\begin{eqnarray*}
 & & L_{0,t} = L_{1,t}=L_{2,t-2}= \cdots = L_{t,0}:=\Lambda, 
\end{eqnarray*}
 where $\Lambda$ is a positive integer, then $(X, {\cal B})$
 forms a simple $t$-design with parameters
 $t-(v,k, \Lambda)$.

 We record the result of the construction discussed above in the 
 following theorem.

\begin{thm}\label{main}
 Let v, k, t be integers with $v > k > t \geq 2$. Let $X$ be a $v$-set
 and let $X=X_1 \cup X_2$ be a partition of $X$ with $|X_1|=v_1$ and
 $|X_2|=v_2$. 
  Let $D_i =(X_1, {\cal B}^{(i)})$ be the 
  complete $i-(v_1,i,1)$ design for $i=0, \ldots, t$  
  and let 
  $D_i =(X_1, {\cal B}^{(i)})$ be a simple $t-(v_1,i, \lambda^{(i)}_t)$
  design for $i=t+1, \ldots, k$.
  Similarly, let 
  $\bar{D}_i =(X_2, \bar{\cal B}^{(i)})$ 
  be the  complete $i-(v_2,i,1)$ design for $i=0, \ldots, t$, and let 
  $\bar{D}_i =(X_2, \bar{\cal B}^{(i)})$ be a 
  simple $t-(v_2,i, \bar{\lambda}^{(i)}_t)$ design for $i=t+1, \ldots, k$.
  Define 
 $${\cal B} = {\cal B}_{(0,k)}\times [u_0] \cup 
             {\cal B}_{(1,k-1)}\times [u_1] 
             \cup   \cdots \cup
             {\cal B}_{(k-1,1)} \times [u_{k-1}] 
             \cup {\cal B}_{(k,0)}\times [u_k], $$   
  where  
$${\cal B}_{(i,{k-i})}= \{ B= B_i \cup \bar{B}_{k-i} \; / \;
   B_i \in {\cal B}^{(i)}, \bar{B}_{k-i} \in \bar{\cal B}^{(k-i)} \}. $$ 
  Assume that 
\begin{eqnarray}\label{Equalities}
  & & L_{0,t}= L_{1,t-1}=L_{2,t-2}= \cdots = L_{t,0}:=\Lambda, 
\end{eqnarray}
  for a positive integer $\Lambda$, where  
  \begin{eqnarray}\label{L-values}
  L_{s,t-s} & = & \sum_{i=0}^k u_i.\lambda^{(i)}_s. 
    \bar{\lambda}^{(k-i)}_{t-s},
  \end{eqnarray}
  $s=0, \ldots , t$, and $u_i\in\{0,1\}$, for $i=0, \ldots, k.$ 
 Then $(X, {\cal B})$ is a simple $t-(v,k, \Lambda)$ design. 
\end{thm}

Two remarks should be included. Firstly, Eq.(\ref{Equalities}) 
always has at least one solution giving rise to the complete 
$t-(v,k, {v-t \choose k-t})$ design.  
In other words, if each ingredient design
is a complete design with its corresponding parameters, then we obtain the
complete design as a result. Secondly, we mainly focus on simple
designs, so we have formulated Theorem \ref{main} accordingly.
But, the construction by no means restricts to simple
$t-$designs. It works for both simple and non-simple designs.
In fact, the construction only uses the ``balance property'' which
depends on the indices $\lambda^{(i)}_t$, and
not on any ``structural property'' of the ingredient designs. 
Thus, if any of the ingredient designs is non-simple, then so is the 
resulting design constructed from a solution of Eq.(\ref{Equalities}).

\section{Applications}
 In this section we illustrate the construction 
 in Theorem \ref{main} through a number of 
 examples which also prove the strength of the method.
 In fact, for some given parameters with $t=4,5,6$, 
 we have constructed a large number of new simple designs.

   In the following we will employ the notation from
   Chapter 4 : $t$-Designs with $t \geq 3$ 
   of the Handbook of Combinatorial Designs.  
   The parameter set $t-(v,k, \lambda)$ of a design will be written as
   $t-(v,k, m\lambda_{\mbox{min}}).$ 
   Since the supplement of a simple
   $t-(v,k, \lambda)$ design is a $t-(v,k, \lambda_{\mbox{max}}-\lambda)$
   design, we usually consider simple $t-(v,k, \lambda)$ designs with
   $\lambda \leq \lambda_{\mbox{max}}/2$.  Thus, the upper limit
   of $m$ of a constructed design will be 
$\mbox{LIM}= \lfloor \lambda_{\mbox{max}}/(2\lambda_{\mbox{min}}) \rfloor.$ 
But, it should be remarked that, when an ingredient design with 
   index $\lambda$ is used, then $\lambda$ can take on
   all possible values, i.e. 
   $\lambda_{\mbox{min}} \leq  \lambda \leq \lambda_{\mbox{max}}$.
  
\subsection{ Simple $5-(36, k, \Lambda)$ designs}

 A detailed example will illustrate the construction.

\subsubsection{Simple $5-(36,10, \Lambda)$ designs} 
 Let $ X= X_1 \cup X_2$ be a partition of the point set $X$ with 
 $|X|=36$ into two subsets $X_1$ and $X_2$ with $|X_1|=|X_2|=18.$ 
 For $i=0,1,2,3,4,5$ let $D_i=(X_1, {\cal B}^{(i)})$ be the complete
 $i-(18,i,1)$ designs. For $i=6,7,8,9,10$ let 
 $D_i=(X_1, {\cal B}^{(i)})$ be a simple $5-(18,i, \lambda^{(i)}_5)$
 design. These designs have the following parameters.   
   \begin{itemize}
    \item[$\bullet$]
     $5-(18,6,\lambda^{(6)}_5) = 5-(18,6,m)$, $m=1,2,\ldots, 13$. 
    \item[$\bullet$]
    $5-(18,7,\lambda^{(7)}_5) =  5-(18,7,m6)$, $m=1,2,\ldots, 13$ 
    \item[$\bullet$]
    $5-(18,8,\lambda^{(8)}_5) =  5-(18,8,m2)$, $m=1,2,\ldots, 143$
    \item[$\bullet$]
    $5-(18,9,\lambda^{(9)}_5) =  5-(18,9,m5)$, $m=1,2,\ldots, 143$
    \item[$\bullet$]
    $5-(18,10,\lambda^{(10)}_5) =  5-(18,10,m9)$, $m=1,2,\ldots, 143$ 
                    (the complement of a $5-(18,8,m2)$).
   \end{itemize}
  Correspondingly, let
  ${\bar D_i}=(X_2, {\bar{\cal B}}^{(i)})$ be simple
  designs defined on $X_2$. 
  We first compute $L_{0,5}$, $L_{1,4}$, $L_{2,3}$.
  We have
  \begin{eqnarray}\label{Eq1} 
  & & L_{s,5-s} = \sum_{i=0}^{10} u_i.\lambda^{(i)}_s. 
    \bar{\lambda}^{(10-i)}_{5-s},
  \end{eqnarray}
  $s=0, \ldots , 5$, and $u_i\in\{0,1\}$ for $i=0, \ldots, 10$.

  Since
  $\bar{\lambda}^{(4)}_5=\bar{\lambda}^{(3)}_5=
   \bar{\lambda}^{(2)}_5=\bar{\lambda}^{(1)}_5=
   \bar{\lambda}^{(0)}_5=0$  and $\bar{\lambda}^{(5)}_5=1,$
   we have  
  \begin{eqnarray*}
   L_{0,5} &=& u_0\lambda^{(0)}_0 \bar{\lambda}^{(10)}_5 +
               u_1\lambda^{(1)}_0 \bar{\lambda}^{(9)}_5 +
               u_2\lambda^{(2)}_0 \bar{\lambda}^{(8)}_5 +
               u_3\lambda^{(3)}_0 \bar{\lambda}^{(7)}_5 +
               u_4\lambda^{(4)}_0 \bar{\lambda}^{(6)}_5 + 
               u_5\lambda^{(5)}_0 \bar{\lambda}^{(5)}_5 \\
           &=& u_0\bar{\lambda}^{(10)}_5 + u_1 18\bar{\lambda}^{(9)}_5+
               u_2 153\bar{\lambda}^{(8)}_5 + u_3 816\bar{\lambda}^{(7)}_5
              + u_4 3060\bar{\lambda}^{(6)}_5 + u_5 8568.               
  \end{eqnarray*}
  Since $\bar{\lambda}^{(3)}_4=\bar{\lambda}^{(2)}_4=
        \bar{\lambda}^{(1)}_4=\bar{\lambda}^{(0)}_4=0$ and 
      ${\lambda}^{(0)}_1=0$, we have
   \begin{eqnarray*}
     L_{1,4} &=& u_1\lambda^{(1)}_1 \bar{\lambda}^{(9)}_4 +
               u_2\lambda^{(2)}_1 \bar{\lambda}^{(8)}_4 +
               u_3\lambda^{(3)}_1 \bar{\lambda}^{(7)}_4 +
               u_4\lambda^{(4)}_1 \bar{\lambda}^{(6)}_4 +
               u_5\lambda^{(5)}_1 \bar{\lambda}^{(5)}_4 +
               u_6\lambda^{(6)}_1 \bar{\lambda}^{(4)}_4 \\
           &=&  u_1 \frac{14}{5}\bar{\lambda}^{(9)}_5+
               u_2 \frac{17\times7}{2}\bar{\lambda}^{(8)}_5 + 
               u_3 \frac{136\times 14}{3}\bar{\lambda}^{(7)}_5 +
               u_4 {680\times 7}\bar{\lambda}^{(6)}_5 + \\ 
             & &   u_5 2380\times 14 +
               u_6 476 {\lambda}^{(6)}_5.
  \end{eqnarray*}  
  Further, since 
  $\bar{\lambda}^{(2)}_3=\bar{\lambda}^{(1)}_3=
   \bar{\lambda}^{(0)}_3={\lambda}^{(0)}_2=
   {\lambda}^{(1)}_2=0$,
   we have  
   \begin{eqnarray*}
     L_{2,3} &=&  u_2\lambda^{(2)}_2 \bar{\lambda}^{(8)}_3 +
               u_3\lambda^{(3)}_2 \bar{\lambda}^{(7)}_3 +
               u_4\lambda^{(4)}_2 \bar{\lambda}^{(6)}_3 +
               u_5\lambda^{(5)}_2 \bar{\lambda}^{(5)}_3 +
               u_6\lambda^{(6)}_2 \bar{\lambda}^{(4)}_3 +
               u_7\lambda^{(7)}_2 \bar{\lambda}^{(3)}_3 \\
           &=& u_2{21 \over 2}\bar{\lambda}^{(8)}_5 + 
               u_3 {16\times 35 \over 2}\bar{\lambda}^{(7)}_5 +
               u_4 {120\times 35}\bar{\lambda}^{(6)}_5 + 
               u_5 560\times 105 + \\
           & &    u_6 {140 \times 15} {\lambda}^{(6)}_5 +
               u_7 56 {\lambda}^{(7)}_5.
  \end{eqnarray*} 
Similarly, we compute
   \begin{eqnarray*}
     L_{3,2} &=&  u_3\lambda^{(3)}_3 \bar{\lambda}^{(7)}_2 +
               u_4\lambda^{(4)}_3 \bar{\lambda}^{(6)}_2 +
               u_5\lambda^{(5)}_3 \bar{\lambda}^{(5)}_2 +
               u_6\lambda^{(6)}_3 \bar{\lambda}^{(4)}_2 +
               u_7\lambda^{(7)}_3 \bar{\lambda}^{(3)}_2 + 
               u_8\lambda^{(8)}_3 \bar{\lambda}^{(2)}_2 \\
           &=& u_3 56\bar{\lambda}^{(7)}_5 +
               u_4 {15\times 140}\bar{\lambda}^{(6)}_5 + 
               u_5 105\times 560 + 
               u_6 {35 \times 120} {\lambda}^{(6)}_5 + \\
            & &   u_7 {35 \times 16 \over 2} {\lambda}^{(7)}_5
                 + u_8 {21 \over 2} {\lambda}^{(8)}_5.
  \end{eqnarray*} 
   \begin{eqnarray*}
   L_{4,1} &=& u_4\lambda^{(4)}_4 \bar{\lambda}^{(6)}_1 +
               u_5\lambda^{(5)}_4 \bar{\lambda}^{(5)}_1 +
               u_6\lambda^{(6)}_4 \bar{\lambda}^{(4)}_1 +
               u_7\lambda^{(7)}_4 \bar{\lambda}^{(3)}_1 +
               u_8\lambda^{(8)}_4 \bar{\lambda}^{(2)}_1 +
               u_9\lambda^{(9)}_4 \bar{\lambda}^{(1)}_1 \\
          &=&    u_4 476\bar{\lambda}^{(6)}_5+
                 u_5 {14\times 2380} + 
                 u_6 {7\times 680}\lambda^{(6)}_5 +
                 u_7 {14\times 136 \over 3}\lambda^{(7)}_5 + \\ 
          & &    u_8 {7 \times 17 \over 2}\lambda^{(8)}_5 + 
                 u_9 {14 \over 5} {\lambda}^{(9)}_5.
  \end{eqnarray*}   
  \begin{eqnarray*}
   L_{5,0} &=& u_5\lambda^{(5)}_5 \bar{\lambda}^{(5)}_0 +
               u_6\lambda^{(6)}_5 \bar{\lambda}^{(4)}_0 +
               u_7\lambda^{(7)}_5 \bar{\lambda}^{(3)}_0 +
               u_8\lambda^{(8)}_5 \bar{\lambda}^{(2)}_0 +
               u_9\lambda^{(9)}_5 \bar{\lambda}^{(1)}_0 + 
               u_{10}\lambda^{(10)}_5 \bar{\lambda}^{(0)}_0 \\
           &=& u_5 8568 + u_6 3060\lambda^{(6)}_5+
               u_7 816\lambda^{(7)}_5 + u_8 153\lambda^{(8)}_5
              + u_9 18\lambda^{(9)}_5 + u_{10}\lambda^{(10)}_5.               
  \end{eqnarray*}
 
   Each set of values  of $u_i \in \{0,1\}$, 
   $i=0, \ldots, 10,$ and 
   $\lambda^{(j)}_5$ and 
   $\bar{\lambda}^{(j)}_5$, $j =6, \ldots, 10,$
   for which the condition
   \begin{eqnarray}
   & & L_{0,5}=L_{1,4}=L_{2,3}=L_{3,2}=L_{4,1}=L_{5,0}:=\Lambda 
   \end{eqnarray}
   is fullfilled for a positive integer $\Lambda$ will 
   yield a simple $5-(36,10, \Lambda)$ design.
   
   Note that a $5-(36,10,\lambda)$ design will be
   written as $5-(36,10, m63)$ with $\lambda_{\mbox{min}}=63$
   and $\lambda_{\mbox{max}}={31 \choose 5}=2697$. So,
   $\mbox{LIM} = \lfloor 2697/{2*63} \rfloor=1348$.  
   By solving Eq.(\ref{Equalities}) above,
   we obtain designs
   for all $m63 \leq 2697.$ 
   Altogether 75 values for $m$ have been found,
   of which 37 values of $m\leq \mbox{LIM}.$
   However, since not
   all  simple $5-(18,i, \lambda^{(i)}_5)$ designs are known to exist,
   for example, $5-(18,6,m)$ designs are known for $m=4,5,6,7,8,9,13$
   only (here $5-(18,6,13)$ is the complete design), we just obtain 
   the following 10 new non-trivial simple $5-(36,10, m63)$ designs
   for $m= 542, 621, 645, 669, 748, 772, 932, 956, 1304, 1328.$ 
   More precisely, Table 1 below shows the details of these 10 solutions.

  \begin{center}
   \begin{tabular}{|r||rrrrrr|} \hline
    $m$ & $\lambda^{(5)}_5$ & $\lambda^{(6)}_5$ & $\lambda^{(7)}_5$ & 
       $\lambda^{(8)}_5$ & $\lambda^{(9)}_5$ & $\lambda^{(10)}_5$ \\ 
        \hline \hline  
    542  &  0  &   5    &   6   &   60   &   210   &   990 \\
    621  &  0  &   6    &   0   &   126  &    75   &   135 \\
    645  &  0  &   6    &   6   &    78  &   275   &   495 \\
    669  &  0  &   6    &  12   &    30  &   475   &   855 \\
    748  &  0  &   7    &   6   &    96  &   340   &     0 \\
    772  &  0  &   7    &  12   &    48  &   540   &   360 \\
    932  &  0  &   9    &   0   &   192  &    60   &   720 \\
    956  &  0  &   9    &   6   &   144  &   260   &  1080 \\
   1304  &  1  &   0    &  66   &   112  &   100   &   792 \\
   1328  &  1  &   0    &  72   &    64  &   300   &  1152 \\ \hline
   \end{tabular}
  \end{center}
    
   An entry $0$ in a column of the table
   implies that $u_i=0$, otherwise $u_i=1$. No values
   for $\bar{\lambda}^{(j)}_5$ are given in the table, because
   we have $\lambda^{(j)}_5 = \bar{\lambda}^{(j)}_5$, $j=6,7,8,9,10,$
   for all these solutions.

\begin{rem}\label{rem1}
{\rm
   In order to simplify the expressions $L_{s,5-s}$ 
   we may introduce the following variables      
   $x_j= u_j\lambda^{(j)}_5$  and $y_j= u_{k-j}{\bar\lambda}^{(j)}_5$
   for $j=6,7,8,9,10$. More precisely,
   
   \[ x_j= \left\{ \begin{array}{ll} 
             0 & \mbox{if $u_j=0$} \\
             \lambda^{(j)}_5 & \mbox{if $u_j=1$}
              \end{array}
      \right. \]
 and
   \[ y_j= \left\{ \begin{array}{ll}
             0 & \mbox{if $u_{k-j}=0$} \\
             \bar{\lambda}^{(j)}_5 & \mbox{if $u_{k-j}=1$}
              \end{array}
      \right. \]
  Thus $L_{s,5-s}$ have much simpler forms, in which $x_j$ and $y_j$
  are allowed to take on the value of zero. For example,
 \begin{eqnarray*} 
    L_{2,3} & = & {21 \over 2}y_8   + 
                  {16\times 35 \over 2}y_7  +
                  {120\times 35}y_6   + 
                  u_5 560\times 105 + 
                  {140 \times 15}x_6  +  
                   56 x_7. \\
    L_{1,4}  & = &  \frac{14}{5} y_9  +
                 \frac{17\times7}{2}y_8   + 
                 \frac{136\times 14}{3}y_7  +
                  {680\times 7}y_6  +  
                   u_5 2380\times 14 +
                   476 x_6  . \\
   L_{0,5}  & = &   y_{10}  + 
                   18 y_9  +
                   153 y_8  + 
                   816 y_7  +
                  3060 y_6 + 
                 u_5 8568.               
\end{eqnarray*} 

}
\end{rem}

\subsubsection{ Simple $5-(36,k,\lambda)$ designs with $11 \leq k \leq 15$}   
 We give a summary of the results from the construction
 of Theorem \ref{main} for simple
 $5-(36,k,\lambda)$ designs for $k=11, \ldots, 15$, for which
 $v_1=v_2=18.$  

 When $v_1=v_2$, we observe that most of the solutions 
 of Eq.(\ref{Equalities}) have the property that
 $\lambda^{(k)}_5=\bar{\lambda}^{(k)}_5$, which we call {\it
 symmetric property}. 
 Thus, assuming symmetric property for solutions of Eq.(\ref{Equalities}) 
 appears to be reasonable. On the other hand, it will reduce
 the search time for solutions enormously.
 For $k=12,13,14,15$ we assume the symmetric property, but even so
 a great number of new designs have been constructed.  

 \begin{itemize}
  
  \item 
        Simple $5-(36,11, \lambda)= 5-(36,11, m21)$ designs with
        $\mbox{LIM}= 17530.$  The construction yields 400 values
        for $m$ with $m\leq \mbox{LIM}$ as solutions for
        Eq.(\ref{Equalities}). The 73 values for $m$ below
  \begin{eqnarray*}
   m &=& 11832, 8712, 8736, 9404, 9416, 9440, 10084, 10120, 10752, 10889,\\
     & & 10913, 11432, 11444, 11456, 11545, 12124, 12136, 12225, 12249, \\
     & & 12261, 12840, 12905, 12929, 12941, 12953, 13496, 14265, 14301, \\
     & & 10676, 10717, 11356, 11397, 12077, 12101, 12781, 12805, 12894, \\
     & & 13396, 13485, 13509, 13574, 14076, 14117, 14189, 14254, 14797, \\   
     & & 14821, 15501, 15614, 16205, 16294, 16861, 16909, 13426, 13450, \\
     & & 14130, 14154, 14834, 14858, 15466, 15538, 16146, 16170, 16271, \\  
     & & 16850, 16874, 16951, 15390, 16070, 16803, 16875, 17483, 17507. \\
  \end{eqnarray*}
   show the constructed simple $5-(36,11, m21)$ designs. 
   Of which 72 values of $m$ yield new designs, except one, $m=13485$, 
    which has been known already. 

   \item
   The results for $k=12,13,14,15 $
   are recorded in the following Table 2.

\begin{center}

\begin{tabular}{|c|c|c|c|} \hline
 Parameters & LIM &\# solutions of Eq.(\ref{Equalities})
            & \# constructed designs  \\ \hline 
$5-(36,12, m15)$  & 87652   &   3261  &  240  \\ 
$5-(36,13, m585)$ & 6742    &   2427  &  359  \\ 
$5-(36,14, m65)$  & 155077  &  26609  &  1926  \\
$5-(36,15, m143)$ & 155077  &  48852  &  4452  \\ \hline           
\end{tabular} 

\end{center}
In Table 2 the figures in column ``\# solutions of 
Eq.(\ref{Equalities})''
are the number of solutions of Eq.(\ref{Equalities}) having 
the symmetric property, whereas those
in column  ``\# constructed designs'' are 
the number of constructed simple designs with
parameters in the first column for $m \leq \mbox{LIM}$. 
The constructed 5-designs are 
derived from solutions of Eq.(\ref{Equalities}) and from known 
simple 5-designs on
18 points as given in \cite{handbook07}.  

\end{itemize}

\begin{rem}
{\rm
We have also applied our method
to constructing $5-(36,k, \Lambda)$ designs for $k=16,17,18$.
In each of these cases we can always construct new designs. 
}
\end{rem}

\vspace{2mm}

\begin{exs} 
 {\rm
 We display some new simple
 5-designs for $k=11,12,13,14,15$ explicitly. All but one design have 
 the symmetric property. The missing values for $\lambda^{(i)}_5$
 and $\bar{\lambda}^{(i)}_5$ in the following examples imply that
 the corresponding designs are not used in the construction.
 Here are the designs.
 \begin{itemize}
 \item $5-(36,11, 11832\times 21)$ with
 $\lambda^{(7)}_5=54$, 
 $\lambda^{(8)}_5=16$,
 $\lambda^{(9)}_5=240$, 
 $\lambda^{(10)}_5=1224$,
       $\bar{\lambda}^{(6)}_5=8$, 
       $\bar{\lambda}^{(7)}_5=12$,
       $\bar{\lambda}^{(8)}_5=108$, 
       $\bar{\lambda}^{(9)}_5=360$. This solution
 does not have the symmetric property.
 
  $5-(36,11, 8712\times 21)$ with
 $\lambda^{(6)}_5= 4$, 
 $\lambda^{(7)}_5= 6$,
 $\lambda^{(8)}_5= 142$,
 $\lambda^{(9)}_5=  40$,
 $\lambda^{(10)}_5= 72$,
 $\lambda^{(11)}_5= 1320,$
and $\bar{\lambda}^{(i)}_5=\lambda^{(i)}_5$, $i=6,7,8,9,10,11.$

\item $5-(36,12, 15337\times 15)$ with
 $\lambda^{(6)}_5= 4$,
 $\lambda^{(7)}_5= 6$,
 $\lambda^{(8)}_5= 30$,
 $\lambda^{(9)}_5= 55$,
 $\lambda^{(10)}_5= 27$,
 $\lambda^{(11)}_5= 660,$
 $\lambda^{(12)}_5= 660,$ 
and $\bar{\lambda}^{(i)}_5=\lambda^{(i)}_5$, $i=6,7,8,9,10,11,12.$

 $5-(36,12, 50490\times 15)$ with
 $\lambda^{(7)}_5= 42$,
 $\lambda^{(8)}_5= 46$,
 $\lambda^{(9)}_5= 135$,
 $\lambda^{(10)}_5= 864,$
 and $\bar{\lambda}^{(i)}_5=\lambda^{(i)}_5$, $i=7,8,9,10.$

\item $5-(36,13, 1347\times 585)$ with
 $\lambda^{(6)}_5= 4$,
 $\lambda^{(7)}_5= 18$,
 $\lambda^{(8)}_5= 48$,
 $\lambda^{(9)}_5= 40$,
 $\lambda^{(10)}_5= 27$,
 $\lambda^{(11)}_5= 396,$
 $\lambda^{(12)}_5= 1716,$ 
 $\lambda^{(13)}_5= 1287,$
and $\bar{\lambda}^{(i)}_5=\lambda^{(i)}_5$, $i=6,7,8,9,10,11,12,13.$

 $5-(36,13, 2448\times 585)$ with
 $\lambda^{(6)}_5= 4$,
 $\lambda^{(7)}_5= 48$,
 $\lambda^{(8)}_5= 48$,
 $\lambda^{(9)}_5= 120$,
 $\lambda^{(10)}_5= 360,$
 and $\bar{\lambda}^{(i)}_5=\lambda^{(i)}_5$, $i=6,7,8,9,10.$

\item $5-(36,14, 20400\times 65)$ with
 $\bar{\lambda}^{(6)}_5= 4$,
 $\bar{\lambda}^{(7)}_5= 30$,
 $\bar{\lambda}^{(9)}_5= 60$,
 $\bar{\lambda}^{(10)}_5= 144$,
and $\lambda^{(i)}_5= \bar{\lambda}^{(i)}_5$, $i=6,7,9,10.$

 $5-(36,14, 19992\times 65)$ with
 $\lambda^{(6)}_5= 4$,
 $\lambda^{(8)}_5= 98$,
 $\lambda^{(9)}_5= 60$,
 $\lambda^{(12)}_5= 1056$,
 and $\bar{\lambda}^{(i)}_5=\lambda^{(i)}_5$, $i=6,8,9,12.$

\item $5-(36,15, 19040\times 143)$ with
 $\lambda^{(6)}_5= 4$,
 $\lambda^{(7)}_5= 6$,
 $\lambda^{(8)}_5= 112$,
 $\lambda^{(9)}_5= 320$,
 $\lambda^{(12)}_5= 528,$ 
and $\bar{\lambda}^{(i)}_5=\lambda^{(i)}_5$, $i=6,7,8,9,12.$

 $5-(36,15, 119952\times 143)$ with
 $\lambda^{(7)}_5= 42$,
 $\lambda^{(8)}_5= 280$,
 $\lambda^{(10)}_5= 1152$,
 $\lambda^{(12)}_5= 792$,
and $\bar{\lambda}^{(i)}_5=\lambda^{(i)}_5$, $i=7,8,10,12.$

 \end{itemize}

}
\end{exs} 

\begin{rem}
{\rm  It is worth mentioning that there may exist different solutions
 to Eq.(\ref{Equalities}) leading to the same value $\Lambda$ for 
 constructed designs. For instance, the following two distinct solutions 
 (a) and (b) of Eq.(\ref{Equalities}) 
 for $t=5, v=36, k=13$:
\begin{itemize}
\item[(a)]
$\lambda^{(6)}_5= 4, \; 
 \lambda^{(7)}_5= 54, \;
 \lambda^{(8)}_5= 128, \; 
 \lambda^{(10)}_5= 729, \;
 \lambda^{(11)}_5= 264, \;$ 
  $\bar{\lambda}^{(i)}_5= \lambda^{(i)}_5, \; i=6,7,8,10,11, $
\item[(b)]
$\lambda^{(6)}_5= 7, \; 
 \lambda^{(7)}_5= 42, \;
 \lambda^{(8)}_5= 64, \; 
 \lambda^{(9)}_5= 240, \; 
 \lambda^{(10)}_5= 288, \;
 \lambda^{(11)}_5= 528, \;$  
   $\bar{\lambda}^{(i)}_5=\lambda^{(i)}_5, \; i=6,7,8,9,10,11,$
\end{itemize}
 lead to simple designs with 
 the same parameters $5-(36,13, 3672 \times 585)$. However,
 they are not isomorphic.

}
\end{rem}
\subsection{Simple $4-(35,k,\Lambda)$ designs with $k=8,9,10$ }
 We shall choose $v_1=17$ and $v_2=18$. 
\subsubsection {$k=8$}
 There is a unique non-trivial solution for Eq.(\ref{Equalities})
 with 
 $\lambda^{(5)}_4=13$, 
 $\lambda^{(7)}_4=264$, 
 $\lambda^{(8)}_4=320$,
     $\bar{\lambda}^{(5)}_4=14$,
     $\bar{\lambda}^{(7)}_4=336$, 
     $\bar{\lambda}^{(8)}_4=448$,
which yields a simple $4-(35,8,448\times 35)$ design.

\subsubsection {$k=9$}
 There are in total 700 non-trivial solutions for Eq.(\ref{Equalities}),
 of which we can construct 452 simple $4-(35,9, \Lambda)$ designs.
 Here are two examples.

(a) $\lambda^{(6)}_4=18$,
    $\lambda^{(7)}_4=38$,
    $\lambda^{(8)}_4=15$,
    $\lambda^{(9)}_4=27$,
         $\bar{\lambda}^{(5)}_4=4$,
         $\bar{\lambda}^{(7)}_4=84$,        
         $\bar{\lambda}^{(8)}_4=133$,
         $\bar{\lambda}^{(9)}_4=42$,
which yields a simple $4-(35,9,369\times 63)$ design.

(b) $\lambda^{(5)}_4=4$,
    $\lambda^{(7)}_4=84$,
    $\lambda^{(8)}_4=50$,
    $\lambda^{(9)}_4=90$,
    $\bar{\lambda}^{(6)}_4=28$,
    $\bar{\lambda}^{(8)}_4=294$,
    $\bar{\lambda}^{(9)}_4=140$,
which yields a simple $4-(35,9,414\times 63)$ design.

\subsubsection {$k=10$}
 There is a huge number of non-trivial solutions 
 for Eq.(\ref{Equalities}) in this case. For instance, with the restriction
 that $\lambda^{(5)}_4=3$, we already have
 constructed 43225 simple $4-(35,10,\Lambda)$ designs (many
 designs have equal value $\Lambda$, but they are not isomorphic).
 Here is an example.

    $\lambda^{(5)}_4=3$,
    $\lambda^{(6)}_4=12$,
    $\lambda^{(7)}_4=6$,
    $\lambda^{(8)}_4=85$,
    $\lambda^{(9)}_4=153$,
    $\lambda^{(10)}_4=612$, 
        $\bar{\lambda}^{(5)}_4=2$,
        $\bar{\lambda}^{(6)}_4=11$,
        $\bar{\lambda}^{(7)}_4=28$, 
        $\bar{\lambda}^{(8)}_4=70$,
        $\bar{\lambda}^{(9)}_4=238$,
        $\bar{\lambda}^{(10)}_4=357$,
which yields a simple $4-(35,10,3043\times 21)$ design.

\subsection{Some simple $6-(46,k,\Lambda)$ designs with $k=13,15$ }
 Some further examples for $6-(46,13,\Lambda)$ and 
 $6-(46,15,\Lambda)$ designs are given here. In
 both cases the ingredient designs are on 23 points, i.e. $v_1=v_2=23.$

 \begin{itemize}
 \item $6-(46,13, 3515\times 1560)$ with
 $\lambda^{(7)}_6=5$, 
 $\lambda^{(8)}_6=40$,
 $\lambda^{(9)}_6=200$, 
 $\lambda^{(10)}_6=700$,
 $\lambda^{(11)}_6=1820$,
 $\lambda^{(12)}_6=3640$, 
 $\lambda^{(13)}_6=5720$,
and $\bar{\lambda}^{(i)}_6=\lambda^{(i)}_6$, $i=7,8,9,10,11,12,13.$
 
 $6-(46,13, 4218\times 1560)$ with
 $\lambda^{(7)}_6=6$,  
 $\lambda^{(8)}_6=48$,
 $\lambda^{(9)}_6=240$, 
 $\lambda^{(10)}_6=840$,
 $\lambda^{(11)}_6=2184$,
 $\lambda^{(12)}_6=4368$,
 $\lambda^{(13)}_6=6864$,
and $\bar{\lambda}^{(i)}_6=\lambda^{(i)}_6$, $i=7,\ldots,13.$

 \item $6-(46,15, 28120\times 2860)$ with
 $\lambda^{(7)}_6=5$, 
 $\lambda^{(8)}_6=136$,
 $\lambda^{(9)}_6=200$, 
 $\lambda^{(10)}_6=700$,
 $\lambda^{(11)}_6=1820$,
 $\lambda^{(12)}_6=3640$, 
 $\lambda^{(13)}_6=5720$,
 $\lambda^{(14)}_6=7150$,
 $\lambda^{(15)}_6=7150$,
and $\bar{\lambda}^{(i)}_6= \lambda^{(i)}_6$, $i=7,\ldots,15.$

\end{itemize}

\vspace{2mm}
\begin{rem}
{\rm
For the cases $ t+1 \leq k \leq 2t-1$ we have observed that 
Eq.(\ref{Equalities}) has 
a unique solution leading to a simple design. 
This is exactly the case, when each ingredient design 
is a complete design, and the resulting design is a complete 
design as well. However, 
when we allow a non-simple design as a resulting design, then 
we may have non-trivial solutions.
} 
\end{rem}

\section{Conclusion}
We have presented a new recursive construction for simple $t-$designs
based on a composition of smaller ingredient designs.
The construction leads to find solutions for the indices of the 
ingredient designs that satisfy a certain set of equalities. With
a small number of examples to demonstrate the strength of the method, 
we have constructed a large amount of new $t-$designs, which were unknown
to date. Clearly the method is very fruitful and powerful.   
We could think of a considerable improvement of the Table for simple 
$t-$designs in the Handbook of Combinatorial Designs, 
when we  would apply this method.


\begin{thebibliography}{99}

\bibitem{Ajoodani96}{\sc S.~Ajoodani-Namini,} Extending large sets
of $t-$designs, {\em J. Combin. Theory A} {\bf 76} (1996) 139--144. 


\bibitem{bjl99}{\sc T.~Beth, D.~Jungnickel and H.~Lenz,}
        {\em Design Theory}, 2nd Edition, Cambridge Univ. Press,
        Cambridge (1999). 

\bibitem{betten95}{\sc A.~Betten, A.~Kerber, A.~Kohnert, R.~Laue, 
            and A.~Wassermann,}
 The discovery of simple 7-designs with automorphism group
    $P\Gamma L(2,32)$. In: {\em Applied algebra, algebraic algorithms
     and error-correcting codes} (eds. G.~Cohen, M.~Giusti and T.~Mora).
      Springer, New York (1995) 131--145.

\bibitem{betten98}{\sc A.~Betten, A.~Kerber, R.~Laue, and A.~Wassermann,}
    Simple 8-designs with small parameters, {\em Des. Codes Crypt.}
    {\bf 15} (1998)  5--27.

\bibitem{betten99}{\sc A.~Betten, R.~Laue, and A.~Wassermann,}
    A Steiner 5-design on 36 points, {\em Des. Codes Crypt.}
    {\bf 17} (1999)  181--186.


\bibitem{bier95a}{\sc J.~Bierbrauer,}
 A family of 4-designs with block size 9, {\em Discr. Math.} {\bf 138}
 (1995) 113-117.

\bibitem{bier95b}{\sc J.~Bierbrauer,}
 A family of 4-designs, {\em Graphs Comb.} {\bf 11} (1995) 209--212.

\bibitem{handbook07}
{\sc C.~J.~Colbourn and J.~H.~Dinitz,} Eds.
{\em Handbook of Combinatorial Designs}, 2nd Edition, CRC Press (2007).


\bibitem{Denniston76}{\sc R.~H.~F.~Denniston,}
 Some new 5-designs, {\em Bull. Lond. Math. Soc.} {\bf 8} (1976) 263--267.


\bibitem{driessen78} {\sc L.~H.~M.~E.~Driessen,}
    $t$-designs, $t\geq 3$, Technical Report, Department of Mathematics,
    Technische Hogeschool Eindhoven, The Netherlands, 1978.
  

\bibitem{Khos-Ajoodani91}
{\sc G.~B.~Khosrovshahi and S.~Ajoodani-Namini,}
Combining $t-$designs, {\em J. Combin. Theory Ser. A}, {\bf 58} (1991), 
26--34. 

\bibitem{Jimbo2011}{\sc M.~Jimbo, Y.~Kunihara, R.~Laue, and M.~Sawa,}
Unifying some infinite families of combinatorial 3-designs,
{\em J. Combin. Theory A} {\bf 118} (2011) 1072--1085.


\bibitem{Kramer76}{\sc E.~S.~Kramer, D.~M.~Mesner,}
 $t-$designs on hypergraphs, {\em Discr. Math.} {\bf 15} (1976) 263--296.


\bibitem{Kramer93}{\sc E.~S.~Kramer, S.~S.~Magliveras, and E.~A.~O'Brien,}
 Some new large sets of t-designs, {\em Australas. J. Combin.} {\bf 7}
 (1993) 189--193.

\bibitem{Kreher93}{\sc D.~L.~Kreher,}
 An infinite family of (simple) 6-designs, {\em J. Combin. Des.} {\bf 1}
(1993) 277-280.
 

\bibitem{Mag84}{\sc S.~S.~Magliveras, and D.~M.~Leavitt,}
 Simple 6-(33,8,36)-designs from $P\Gamma L_2(32)$,
 {\em  Computational Group Theory,} Academic Press, New York
 (1984) 337--352.

\bibitem{Mag87}{\sc S.~S.~Magliveras, and T.~E.~Plambeck,} 
 New infinite families of simple 5-designs, {\em J. Combin. Theory A}
 {\bf 44} (1987) 1--5.


\bibitem{Tuan2002}{\sc Ngo Dac Tuan,}  
Simple non-trivial designs with an arbitrary automorphism group,
{\em J. Combin. Theory A} {\bf 100} (2002) 403--408.

\bibitem{Sebille2001}{\sc M.~Sebille,}
There exists a simple non-trivial $t-$design with an arbitrarily large
 automorphism group for every $t$, {\em Des. Codes Crypt.} {\bf 22}
 (2001) 203--206.


\bibitem{Teirlinck87}{\sc L.~Teirlinck,}
 Non-trivial $t-$designs without repeated blocks exist for all $t$,
 {\em Discr. Math.} {\bf 65} (1987) 301--311

\bibitem{Teirlinck89}{\sc L.~Teirlinck,}
Locally trivial $t-$designs and $t-$designs without repeated blocks,
{\em Discr. Math.} {\bf 77} (1989) 345--356.

\bibitem{TvT84}{\sc Tran van Trung,}
 The existence of an infinite family of simple 5-designs,
 {\em Math. Zeitschr.} {\bf 187} (1984) 285--287.

\bibitem{TvT86}{\sc Tran van Trung,}
On the construction of $t-$designs and the existence of some new infinite
families of simple 5-designs, {\em Arch. Math.} {\bf 47} (1986) 187--192.

\bibitem{TvT2001}{\sc Tran van Trung,}
Recursive constructions for 3-designs and resolvable 3-designs,
{\em J. Stat. Plann. Infer.} {\bf 95} (2001) 341--358.

\bibitem{QRWu91}{\sc Qui-rong Wu,}
 A note on extending $t-$designs,
 {\em Australas. J. Combin.} {\bf 4} (1991) 229--235.

\end{thebibliography}
\end{document}